\documentclass[12pt]{amsart} 
\usepackage{amssymb,latexsym}
\usepackage{enumerate}

\makeatletter
\@namedef{subjclassname@2010}{%
  \textup{2010} Mathematics Subject Classification}
\makeatother

\newtheorem{thm}{Theorem}[section]
\newtheorem{cor}[thm]{Corollary}
\newtheorem{lem}[thm]{Lemma}

\newtheorem{prop}[thm]{Proposition}

\theoremstyle{definition}

\newtheorem{rem}[thm]{Remark}

\numberwithin{equation}{section}

\newcommand{\bC}{{\mathbb C}}
\newcommand{\bK}{{\mathbb K}}
\newcommand{\bL}{{\mathbb L}}
\newcommand{\bQ}{{\mathbb Q}}

\newcommand{\bZ}{{\mathbb Z}}

\newcommand{\cA}{{\mathcal A}}
\newcommand{\cC}{{\mathcal C}}
\newcommand{\cD}{{\mathcal D}}
\newcommand{\cN}{{\mathcal N}}

\newcommand{\core}{{\rm core}}
\newcommand{\Gal}{{\rm Gal}}
\newcommand{\id}{{\rm identity}}

\begin{document}

%%%%% To ease editing, for IMPAN journals add:

\baselineskip=17pt

%%%%%%%%%%%

%% In the running head, replace first names by initials 
%% and give an abbreviation of the title.

\title[Effective Irrationality Measures]{Effective Irrationality Measures and Approximation by Algebraic Conjugates}

\author{Paul M. Voutier}
\address{London, UK}
\email{Paul.Voutier@gmail.com}

\date{}

\begin{abstract}
In this paper, we present a result on using algebraic conjugates to form a
sequence of approximations to an algebraic number, and in this way obtain
effective irrationality measures for related algebraic numbers. From this
result, we are able to generalise Thue's Fundamentaltheorem.
\end{abstract}

%    \subjclass is required.

\maketitle

\section{Introduction}

In a recent article \cite{Vout}, we investigated Thue's Fundamentaltheorem
\cite{Thue}, showing when it can be used and how to use it in these cases.
Using the notation of Theorems~1 and 2 of \cite{Vout}, we also showed that the
case when $[\bK(\beta_{1}):\bK]=1$ is equivalent to the ``usual'' hypergeometric
method (see Corollary~1 of \cite{Vout}), where, here and in what follows, $\bK$
is either $\bQ$ or an imaginary quadratic field.

We also considered the case of $[\bK(\beta_{1}):\bK]=2$ in \cite{Vout}.
The approximants $P_{r}(x)$ and $Q_{r}(x)$ that we defined in Lemma~3.3 of
\cite{Vout} have a particularly nice form: an algebraic number plus or minus
its algebraic conjugate. This raises the intriguing question of why.

We address that question here and show that the form of $P_{r}(x)$ and $Q_{r}(x)$
arises from the fact that Thue's Fundamentaltheorem is a special case of the
application to hypergeometric polynomials of a simple observation regarding
diophantine approximations.

We present this observation here along with a generalisation and extension of
Thue's Fundamentaltheorem. In the notation of \cite{Vout}, we are now able to
consider more general expressions in place of $W(x)$ (see also Remark~\ref{rem:w})
as well as more general expressions for the denominator of $\cA(x)$. There are
also smaller improvements such as the consideration of powers $m/n$ rather than
just $1/n$, simplification of the numerator of $\cA(x)$,\ldots

The cost of these improvements is merely in the constant $c$ that appears in
our results below. The irrationality measure, $\kappa$, itself remains unchanged.

\section{Notation}

For positive integers $m$ and $n$ with $0<m<n$, $(m,n)=1$ and a non-negative
integer $r$, we put
\begin{displaymath}
X_{m,n,r}(x) = {} _{2}F_{1}(-r,-r-m/n;1-m/n;x),  
\end{displaymath}
where $_{2}F_{1}$ denotes the classical hypergeometric function. 

We use $X_{m,n,r}^{*}$ to denote the homogeneous polynomials derived from these
polynomials, so that  
\begin{displaymath}
X_{m,n,r}^{*}(x,y) = y^{r} X_{m,n,r}(x/y).  
\end{displaymath}

We let $D_{n,r}$ denote the smallest positive integer such that $D_{n,r} X_{m,n,r}(x)$
has rational integer coefficients.

For a positive integer $d$, we define $N_{d,n,r}$ to be the greatest common
divisor of the numerators of the coefficients of $X_{m,n,r}(1-dx)$.

We will use $v_{p}(x)$ to denote the largest power of a prime $p$ which divides
into the rational number $x$. With this notation, for positive integers $d$ and
$n$, we put
\begin{equation}
\label{eq:Ndn}
\cN_{d,n} =\prod_{p|n} p^{\min(v_{p}(d), v_{p}(n)+1/(p-1))}.
\end{equation}

For any complex number $w$, we can write $w=|w|e^{i \varphi}$, where $|w| \geq 0$
and $-\pi < \varphi \leq \pi$ (with $\varphi=0$, if $w=0$). With such a
representation, unless otherwise stated, $w^{m/n}$ will signify
$\left( |w|^{1/n} \right)^{m} e^{im \varphi/n}$ for positive integers $m$ and
$n$, where $|w|^{1/n}$ is the unique non-negative $n$-th root of $|w|$.

Lastly, following the function name in PARI, we define $\core(n)$ to be the
unique squarefree divisor, $n_{1}$, of $n$ such that $n/n_{1}$ is a perfect square.

\section{Results}

\begin{prop}
\label{prop:1}
Let $\bK$ be either $\bQ$ or an imaginary quadratic field. Let $s \geq 2$ be a
positive integer and $\bL$ be a number field with $[\bL:\bK]=s$.

Let $\theta_{1}=1$, $\theta_{2}$, \ldots, $\theta_{s} \in \bC$ be linearly-independent
over $\bK$ and let $\sigma_{1}=\id$, \ldots, $\sigma_{s}$ be the $s$
embeddings of $\bL$ into $\bC$ that fix $\bK$.

Suppose that there exist real numbers $k_{0},l_{0} > 0$ and $E,Q > 1$ such that
for all non-negative integers $r$, there are algebraic integers $p_{r} \in \bL$
with $\max_{1 \leq i \leq s} |\sigma_{i}(p_{r})|<k_{0}Q^{r}$.

Let $\beta$ and $\gamma$ be algebraic integers in $\bL$.

\noindent
{\rm (i)} Assume that
$\sum_{1 \leq i,j \leq s} \left\{ \sigma_{i}(\beta)\sigma_{j}(\gamma)-\sigma_{j}(\beta)\sigma_{i}(\gamma) \right\}
\sigma_{i}(p_{r})\sigma_{j}(p_{r+1}) \neq 0$ and
$\max_{2 \leq i \leq s} \left| p_{r}\theta_{i}-\sigma_{i}(p_{r}) \right| < l_{0}E^{-r}$.
Put
$$
\alpha = \frac{\sum_{i=1}^{s} \sigma_{i}(\beta)\theta_{i}}{\sum_{i=1}^{s} \sigma_{i}(\gamma)\theta_{i}}.
$$

For any algebraic integers $p$ and $q$ in $\bK$ with $q \neq 0$, we have
$$
\left| \alpha - \frac{p}{q} \right| > \frac{1}{c |q|^{\kappa+1}},
$$
where
$$
c = 2 \left( \sum_{i=1}^{s} \left| \sigma_{i}(\gamma) \right| \right) k_{0}Q
    \max \left\{ E, 2 \left( \sum_{i=2}^{s} \left| \sigma_{i}(\beta) - \alpha \sigma_{i}(\gamma) \right| \right)
    l_{0}E \right\}^{\kappa}
\mbox{ and } \kappa = \frac{\log Q}{\log E}.
$$

\noindent
{\rm (ii)} For $s=2$, assume that $\beta/\gamma, p_{r}/p_{r+1} \not\in \bK$,
and either $\left| p_{r}\theta_{2}-\sigma_{2}(p_{r}) \right| < l_{0}E^{-r}$
or
$\left| -p_{r}\theta_{2}-\sigma_{2}(p_{r}) \right| < l_{0}E^{-r}$. Put
$$
\alpha = \frac{\sigma_{2}(\beta) \theta_{2} \pm \beta}
{\sigma_{2}(\gamma) \theta_{2} \pm \gamma},
$$
where the operation in the numerator matches the operation in the denominator.
If $\bK=\bQ$, then let $\tau=1$, else let $\tau$ be an algebraic integer in
$\bK$ such that $\bL=\bK(\sqrt{\tau})$.

For any algebraic integers $p$ and $q$ in $\bK$ with $q \neq 0$, we have
$$
\left| \alpha - \frac{p}{q} \right| > \frac{1}{c |q|^{\kappa+1}},
$$
where
$$
c = 2 |\sqrt{\tau}| \left( |\gamma| +|\sigma_{2}(\gamma)| \right) k_{0}Q
    \max \left\{ E, 2|\sqrt{\tau}| |\sigma_{2}(\beta)-\alpha \sigma_{2}(\gamma)| l_{0}E \right\}^{\kappa}
\mbox{ and } \kappa = \frac{\log Q}{\log E}.
$$
\end{prop}

We will use part~(ii) of this Proposition to prove the following theorems.

\begin{thm}
\label{thm:general-hypg}
\footnote{Note that our Theorems and Corollary here correct a small error in
Theorems~2.1, 2.4 and Corollary~2.7 of \cite{Vout}, where $\max(1,\ldots$ in
the expressions for $c$ should read $\max(E,\ldots$.}
Let $\bK$ be either $\bQ$ or an imaginary quadratic field. Let $\bL$ be a
number field with $[\bL:\bK]=2$ and let $\sigma$ be the non-trivial element of
$\Gal(\bL/\bK)$. If $\bK=\bQ$, then let $\tau=1$, else let $\tau$ be an
algebraic integer in $\bK$ such that $\bL=\bK(\sqrt{\tau})$. Let $\beta$,
$\gamma$, $\eta$ be algebraic integers in $\bL$.

Let $g$ be an algebraic number such that $\eta/g$ and $\sigma(\eta)/g$ are
algebraic integers $($not necessarily in $\bL)$. For each non-negative integer
$r$, let $h_{r}$ be a non-zero algebraic integer with $h_{r}/g^{r} \in \bK$
and $|h_{r}| \leq h$ for some fixed positive real number $h$. Let $d$ be the
largest positive rational integer such that $(\sigma(\eta)-\eta)/(dg)$
is an algebraic integer and let $\cC_{n}$ and $\cD_{n}$ be positive real
numbers such that
\begin{equation}
\label{eq:num-denom-bnd-a}
\max \left( 1, \frac{\Gamma(1-m/n) \, r!}{\Gamma(r+1-m/n)},
\frac{n\Gamma(r+1+m/n)}{m \Gamma(m/n)r!} \right)
\frac{D_{n,r}}{N_{d,n,r}} < \cC_{n} \left( \frac{\cD_{n}}{\cN_{d,n}} \right)^{r}
\end{equation}
holds for all non-negative integers $r$.

Put
\begin{eqnarray*}
\alpha & = & \frac{\beta(\eta/\sigma(\eta))^{m/n} \pm \sigma(\beta)}
                  {\gamma(\eta/\sigma(\eta))^{m/n} \pm \sigma(\gamma)}, \\
E      & = & \left\{ \frac{\cD_{n}}{|g|\cN_{d,n}} \min \left( \left| \sqrt{\eta}-\sqrt{\sigma(\eta)} \right|^{2},
             \left| \sqrt{\eta}+\sqrt{\sigma(\eta)} \right|^{2} \right) \right\}^{-1}, \\
Q      & = & \frac{\cD_{n}}{|g|\cN_{d,n}} \max \left( \left| \sqrt{\eta}-\sqrt{\sigma(\eta)} \right|^{2},
             \left| \sqrt{\eta}+\sqrt{\sigma(\eta)} \right|^{2} \right), \\
\kappa & = & \frac{\log Q}{\log E} \mbox{ and } \\
c      & = & 4h |\sqrt{\tau}| \left( |\gamma| + |\sigma(\gamma)| \right) \cC_{n}Q \\
       &   & \times \max \left\{ E, 5h|\sqrt{\tau}|\left| 1- (\eta/\sigma(\eta))^{m/n} \right|
                                    |\beta-\alpha \gamma| \cC_{n} E \right\}^{\kappa},
\end{eqnarray*}
where the operation in the numerator of the definition of $\alpha$ matches the
operation in its denominator.

If $E > 1$ and either $0 < \eta/\sigma(\eta) < 1$ or $|\eta/\sigma(\eta)|=1$ with
$\eta/\sigma(\eta) \neq -1$, then 
\begin{equation}
\label{eq:gen-result}
\left| \alpha - p/q \right| > \frac{1}{c |q|^{\kappa+1}} 
\end{equation}
for all algebraic integers $p$ and $q$ in $\bK$ with $q \neq 0$. 
\end{thm}

\begin{rem}
\label{rem:w}
Observe that in our definition of $\alpha$, we take the $n$-th root of
$\eta/\sigma(\eta)$. However, this is more general than it may first
appear. It can be applied to any quantity $\mu \eta/\sigma(\eta)$ where
$\mu \in \bL$ and $\mu=\nu/\sigma(\nu)$ for some $\nu \in \bL$.

For example, although in Thue's Fundamentaltheorem we take the $n$-th root of
$-\eta/\sigma(\eta)$, it, and its generalisations, still follows from our
results. Suppose $\bL=\bK(\sqrt{\tau})$ and put $\eta'=\sqrt{\tau}\eta$, then
$-\eta/\sigma(\eta)=\eta'/\sigma(\eta')$,
so we can express $-\eta/\sigma(\eta)$ in the form here (i.e., take $\mu=-1$
and $\nu=\sqrt{\tau}$ in the above notation). There appears to be an
extra factor of $\sqrt{\tau}$ that will arise in our expressions for $E$ and
$Q$, but these are in fact cancelled out since $g$ also increases by a factor
of $\sqrt{\tau}$, so $\kappa$ is unaffected.

Similarly, if $\bK \neq \bQ(i)$ and $\bL=\bK(i)$. Then
$i\eta/\sigma(\eta)=\eta'/\sigma(\eta')$, where $\eta'=(1+i)\eta$.

Also, if $\bK \neq \bQ(\sqrt{-3})$ and $\bL=\bK(\sqrt{-3})$. Then
$\zeta_{3}\eta/\sigma(\eta)=\eta'/\sigma(\eta')$, where $\eta'=(1-\sqrt{-3})\eta/2$.
And $\zeta_{6}\eta/\sigma(\eta)=\eta'/\sigma(\eta')$, where $\eta'=(3+\sqrt{-3})\eta$.

As for the other roots of unity of degree at most $4$ over $\bQ$, it can be
shown, via algebraic manipulation, that this is not possible for $\zeta_{8}$
and $\zeta_{12}$. And since $\bQ(\zeta_{5})$ contains no subfields besides
$\bQ$ and $\bQ(\sqrt{5})$, we cannot consider $\zeta_{5}\eta/\sigma(\eta)$.
\end{rem}

\begin{rem}
From Lemma~7.4 of \cite{Vout}, the inequality (\ref{eq:num-denom-bnd-a}) holds
for $\cC_{n}$ and $\cD_{n}$ as in \cite{Vout} and hence it does not impose any
constraint.
\end{rem}

\begin{thm}
\label{thm:general-hypg-unitdisk}
Let $\bK$ be an imaginary quadratic field and $\alpha, \beta, \gamma, \eta$,
$\tau$, $\sigma$, $d, g$, $h$, $n$, $\cC_{n}$, $\cD_{n}, \cN_{d,n}$ be as in
Theorem~$\ref{thm:general-hypg}$.

Put
\begin{eqnarray*}
E      & = & \frac{4|g|\cN_{d,n}}{\cD_{n}} \frac{\left( |\eta| - |\sigma(\eta) - \eta | \right)}{|\sigma(\eta) - \eta |^{2}}, \\
Q      & = & \frac{2\cD_{n}}{|g|\cN_{d,n}} \left( \left| \eta \right| + \left| \sigma(\eta) \right| \right), \\
\kappa & = & \frac{\log Q}{\log E} \mbox{ and } \\        
c      & = & 4h |\sqrt{\tau}|\left( |\gamma| + |\sigma(\gamma)| \right) \cC_{n}Q \\
       &   & \times \max \left\{ E, 2h|\sqrt{\tau}|\left| 1- (\eta/\sigma(\eta))^{m/n} \right|
                                    |\beta - \alpha \gamma| \cC_{n} E \right\}^{\kappa}.
\end{eqnarray*}

If $E>1$ and $\max \left( |1-\eta/\sigma(\eta)|, |1-\sigma(\eta)/(\eta)| \right)<1$,
then
\begin{equation}
\label{eq:gen-unitdisk-result}
\left| \alpha - p/q \right| > \frac{1}{c |q|^{\kappa+1}} 
\end{equation}
for all algebraic integers $p$ and $q$ in $\bK$ with $q \neq 0$. 
\end{thm}

\begin{rem}
The condition that $\bK$ be an imaginary quadratic field is no restriction since
the case of $\bK=\bQ$ is completely covered by Theorem~\ref{thm:general-hypg}.
\end{rem}

We now present a corollary of Theorem~\ref{thm:general-hypg} when $\bK=\bQ$.

\begin{cor}
\label{cor:cor-1}
Let $\bK=\bQ$ and $\alpha, \beta, \gamma, \eta$, $\sigma$, $n$, $\cC_{n}$,
$\cD_{n}, \cN_{d,n}$ be as in Theorem~$\ref{thm:general-hypg}$. Suppose that
$\eta=(u_{1} + u_{2} \sqrt{t})/2$ where $t, u_{1}, u_{2} \in \bZ$ and $t \neq 0$.

Put
\allowdisplaybreaks
\begin{eqnarray*}
g_{1}  & = & \gcd \left( u_{1}, u_{2} \right), \\
g_{2}  & = & \gcd(u_{1}/g_{1}, t), \\
g_{3}  & = & \left\{ 	\begin{array}{ll}
             1 & \mbox{if $t \equiv 1 \bmod 4$ and $(u_{1}-u_{2})/g_{1} \equiv 0 \bmod 2$}, \\
             2 & \mbox{if $t \equiv 3 \bmod 4$ and $(u_{1}-u_{2})/g_{1} \equiv 0 \bmod 2$},\\
             4 & \mbox{otherwise,}
             \end{array}
                       \right. \\
g_{4}  & = & \gcd \left( \core(tg_{2}g_{3}), \frac{n}{\gcd((u_{2}/g_{1})\sqrt{tg_{3}/g_{2}/\core(tg_{2}g_{3})}, n)} \right), \\
g_{5}  & = & \left\{ 	\begin{array}{ll}
             2 & \mbox{if $2|n$ and $v_{2} \left( u_{2}^{2}tg_{3}/(g_{1}^{2}g_{2}) \right) = v_{2} \left( 2n^{2} \right)$}, \\
             1 & \mbox{otherwise}
             \end{array}
                       \right.
\hspace{3.0mm} \mbox{and} \\
g      & = & \frac{g_{1}\sqrt{g_{2}}}{\sqrt{g_{3}g_{4}g_{5}}},\\
E      & = & \frac{|g|\cN_{d,n}}{\cD_{n}\min \left( \left| u_{1} \pm \sqrt{u_{1}^{2}-u_{2}^{2}t} \right| \right)}, \\
Q      & = & \frac{\cD_{n}\max \left( \left| u_{1} \pm \sqrt{u_{1}^{2}-u_{2}^{2}t} \right| \right)}{|g|\cN_{d,n}}, \\
\kappa & = & \frac{\log Q}{\log E} \mbox{ and } \\        
c      & = & 4 \sqrt{|2t|} \left( |\gamma| + |\sigma(\gamma)| \right) \cC_{n} Q \\
       &   & \times \left( \max \left( E, 5 \sqrt{|2t|} \left| 1- (\eta/\sigma(\eta))^{m/n} \right|
             |\beta-\alpha\gamma| \cC_{n}E \right) \right)^{\kappa},
\end{eqnarray*}
where $d$ is the largest positive rational integer such that $u_{2}\sqrt{t}/(dg)$
is an algebraic integer.

If $E > 1$ and either $0 < \eta/\sigma(\eta) < 1$ or $|\eta/\sigma(\eta)|=1$
with $\eta/\sigma(\eta) \neq -1$, then
\begin{equation}
\label{eq:cor-1-result}
\left| \alpha - p/q \right| > \frac{1}{c |q|^{\kappa+1}} 
\end{equation}
for all rational integers $p$ and $q$ with $q \neq 0$. 
\end{cor}

\begin{rem}
\label{rem:g}
The factors, $g_{i}$, used to construct $g$ each arise in natural and distinct
ways. $g_{1}$ through $g_{3}$ provide ways to remove common factors from $\eta$
and $\sigma(\eta)$. $g_{4}$ and $g_{5}$ arise from the interplay of $d$ and $g$:
under some circumstances (captured by $g_{4}$ and $g_{5}$), decreasing $g$ can
increase $d$ and hence $\cN_{d,n}$ by more to provide a net benefit.
\end{rem}

\begin{rem}
Using the same argument as in the proof of this Corollary, we can also improve
Corollary~2.7 of \cite{Vout}, replacing $g_{4}$ there by
$$
\gcd \left( \core(g_{2}g_{3}), \frac{n}{\gcd((u_{1}/g_{1})\sqrt{g_{2}/g_{3}/\core(g_{2}g_{3})}, n)} \right)
$$
and adding an appropriate version of the $g_{5}$ above by setting $g_{5}=2$ if
$2|n$ and
$v_{2} \left( u_{1}g_{3}/(g_{1}^{2}g_{2}) \right) = v_{2} \left( 2n^{2} \right)$
and setting $g_{5}=1$ otherwise, since the definition of $d$ in Corollary~2.7
of \cite{Vout} uses $u_{1}/(dg)$ rather than $u_{2}\sqrt{t}/(dg)$ as here.

This improved version of Corollary~2.7 of \cite{Vout} will yield the same
results as in the Corollary here together with Remark~\ref{rem:w}.
\end{rem}

\section{Preliminary Lemmas}

The next lemma contains the relationship that allows the hypergeometic method to
provide good sequences of rational approximations.

\begin{lem}
\label{lem:relation}
For any positive integers $m$ and $n$ with $(m,n)=1$, any non-negative integer
$r$ and for any complex number $z$ that is not a negative number and not zero,
\begin{equation}
\label{eq:approx}
z^{m/n} z^{r} X_{m,n,r}(z^{-1}) - X_{m,n,r}(z) = (z-1)^{2r+1} R_{m,n,r}(z), 
\end{equation}
where
$$
(z-1)^{2r+1} R_{m,n,r}(z) = \frac{\Gamma(r+1+m/n)}{r!\Gamma(m/n)}
\int_{1}^{z} (1-t)^{r}(t-z)^{r}t^{m/n-r-1}dt.
$$
\end{lem}

\begin{rem}
Note that the expression $(z-1)^{2r+1} R_{m,n,r}(z)$ here is the same as the
$R_{m,n,r}(z)$ defined in Lemma~7.1 of \cite{Vout}.
\end{rem}

\begin{proof}
This is shown in the case of $m=1$ in the proof of Lemma~2.3 of \cite{CV}.
The proof for arbitrary $m$ is identical.
\end{proof}

\begin{lem}
\label{lem:approx}
Let $\theta \in \bC$ and let $\bK$ be either $\bQ$ or an imaginary quadratic field.
Suppose that there exist real numbers $k_{0},l_{0} > 0$ and $E,Q > 1$ such that
for all non-negative integers $r$, there are algebraic integers $p_{r}$ and $q_{r}$
in $\bK$ with $|q_{r}|<k_{0}Q^{r}$ and $|q_{r}\theta-p_{r}| \leq l_{0}E^{-r}$
satisfying $p_{r}q_{r+1} \neq p_{r+1}q_{r}$. Then for any algebraic integers $p$
and $q$ in $\bK$ with $q \neq 0$, we have
$$
\left| \theta - \frac{p}{q} \right| > \frac{1}{c |q|^{\kappa+1}},
\mbox{ where } c=2k_{0}Q \left( \max ( 1, 2l_{0}) E \right)^{\kappa}
\mbox{ and } \kappa = \frac{\log Q}{\log E}.
$$

Moreover, if $p/q \neq p_{i}/q_{i}$ for any non-negative integer $i$, then we
can put $c=2k_{0} \left( \max ( 1, 2l_{0}) E \right)^{\kappa}$.
\end{lem}

\begin{proof}
This follows from Lemma~6.1 of \cite{Vout}.

There we proved a similar result for $|q| \geq 1/(2l_{0})$ and
$c=2k_{0}Q(2l_{0}E)^{\kappa}$. Here we merely observe that if we replace $l_{0}$
with $\max ( 0.5, l_{0})$, then all the hypotheses of the Lemma still hold.
Moreover, $1/(2\max ( 0.5, l_{0})) \leq 1$, so the result holds for all non-zero
algebraic integers $q \in \bK$.

The last statement in the Lemma follows since the $Q$ which appears in the
expression for $c$ in the statement of Lemma~6.1 of \cite{Vout} arises only
from consideration of the case $p/q=p_{i}/q_{i}$ for some positive integer $i$. 
\end{proof}

\section{Proof of Proposition~\ref{prop:1}}

Assume that we have a sequence of $p_{r}$'s satisfying the hypotheses of
Proposition~\ref{prop:1}.

\noindent
{\bf (i)} Suppose we have $p_{r}\theta_{i}-\sigma_{i}(p_{r})=\delta_{i,r}$ for
each $i=1$, \ldots, $s$. Then we can write
$$
\alpha = \frac{\sum_{i=1}^{s} \sigma_{i}(\beta)(\delta_{i,r}+\sigma_{i}(p_{r}))}
         {\sum_{i=1}^{s} \sigma_{i}(\gamma)(\delta_{i,r}+\sigma_{i}(p_{r}))}.
$$
and hence
$$
\alpha\sum_{i=1}^{s} \sigma_{i}(\gamma p_{r})
- \sum_{i=1}^{s} \sigma_{i}(\beta p_{r})
= \sum_{i=2}^{s} \left( \sigma_{i}(\beta) - \alpha\sigma_{i}(\gamma) \right)\delta_{i,r},
$$
since $\delta_{1,r}=0$.

Put $p_{r}'=\sum_{i=1}^{s} \sigma_{i}(\beta p_{r})$ and
$q_{r}'=\sum_{i=1}^{s} \sigma_{i}(\gamma p_{r})$. Note that both
$p_{r}'$ and $q_{r}'$ are algebraic integers in $\bK$.

Observe that
$$
\left| \alpha q_{r}' - p_{r}' \right|
< l_{0} \left( \sum_{i=2}^{s} \left| \sigma_{i}(\beta)-\alpha \sigma_{i}(\gamma) \right| \right)
E^{-r}
$$
and
$$
\left| q_{r}' \right|
\leq k_{0} \left( \sum_{i=1}^{s} \left|\sigma_{i}(\gamma) \right| \right)
Q^{r}.
$$

Since
$$
p_{r}'q_{r+1}'-p_{r+1}'q_{r}'
=
\sum_{1 \leq i,j \leq s} \left\{ \sigma_{i}(\beta)\sigma_{j}(\gamma)-\sigma_{j}(\beta)\sigma_{i}(\gamma) \right\}
\sigma_{i}(p_{r})\sigma_{j}(p_{r+1})
\neq 0
$$
by our assumption in the statement of the Proposition, we can apply
Lemma~\ref{lem:approx} with $p_{r}'$ and $q_{r}'$ instead of $p_{r}$ and $q_{r}$,
respectively, to complete the proof in this case.

\noindent
{\bf (ii)} Suppose we have $\zeta_{2}p_{r}\theta_{2}-\sigma_{2}(p_{r})
=\delta_{2,r}$ for some square root of $1$, $\zeta_{2}$, fixed for a given
value of $r$. As above, we can write
$$
\alpha \left\{ \sigma_{2} \left( \gamma p_{r} \right) \pm \zeta_{2}\gamma p_{r} \right\}
- \left\{ \sigma_{2} \left( \beta p_{r} \right) \pm \zeta_{2}\beta p_{r} \right\}
= \delta_{2,r} \left( \sigma_{2} \left(\beta\right) - \alpha \sigma_{2} \left( \gamma \right) \right).
$$

We break the proof into two cases depending on the value of $\zeta_{2}$.

\vspace{1.0mm}

\noindent
\underline{Case 1:} $\pm \zeta_{2}=1$

This case is identical to part~(i) with $s=2$.

Note that in this case ($s=2$), the condition in part~(i) reduces to
$$
\left( \sigma_{2}(\beta)\gamma - \beta\sigma_{2}(\gamma) \right)
\left( \sigma_{2}(p_{r})p_{r+1} - p_{r}\sigma_{2}(p_{r+1}) \right) \neq 0.
$$

This is true under the conditions we have stipulated here, namely
$\beta/\gamma \not\in \bK$ and $p_{r}/p_{r+1} \not\in \bK$
(since the fixed field of $\sigma_{2}$ is $\bK$).

Also since $|\tau| \geq 1$, our definition of $c$ is valid.

\vspace{1.0mm}

\noindent
\underline{Case 2:} $\pm \zeta_{2}=-1$

We break this case into two subcases.

\vspace{1.0mm}

\noindent
\underline{Case 2(i):} $\pm \zeta_{2}=-1$ and $\bK=\bQ$

If $\bK=\bQ$, then we can write $\beta p_{r}=(a+b\sqrt{t})/2$ for some choice
of rational integers $a$, $b$ and $t$ with $t \neq 0$. Hence
$\beta p_{r} - \sigma_{2} \left( \beta p_{r} \right)=b\sqrt{t}$ and
$(\beta p_{r} - \sigma_{2} \left( \beta p_{r} \right))/\sqrt{t} \in \bZ$.
Similarly,
$(\gamma p_{r} - \sigma_{2} \left( \gamma p_{r} \right))/\sqrt{t} \in \bZ$.

In this case, we put
$q_{r}' = \left( \gamma p_{r} - \sigma_{2} \left( \gamma p_{r} \right) \right)/\sqrt{t}$
and $p_{r}' = \left( \beta p_{r} - \sigma_{2} \left( \beta p_{r} \right) \right)/\sqrt{t}$
and observe that
$$
\left| \alpha q_{r}' - p_{r}' \right|
  <  \frac{l_{0} |\sigma_{2}(\beta)-\alpha \sigma_{2}(\gamma)|}{|\sqrt{t}|} E^{-r}
\leq l_{0} |\sqrt{\tau}| |\sigma_{2}(\beta)-\alpha \sigma_{2}(\gamma)| E^{-r}
$$
and
$$
\left| q_{r}' \right|
\leq \frac{k_{0} \left( |\gamma| + |\sigma_{2}(\gamma)| \right)}{|\sqrt{t}|} Q^{r}
\leq k_{0} |\sqrt{\tau}| \left( |\gamma| + |\sigma_{2}(\gamma)| \right) Q^{r},
$$
since $|t| \geq 1$.

\vspace{1.0mm}

\noindent
\underline{Case 2(ii):} $\pm \zeta_{2}=-1$ and $\bK$ is an imaginary quadratic field

If $\bK$ is an imaginary quadratic field, then $\beta p_{r}=a+b\sqrt{\tau}$ for
some $a, b \in \bK$ and where $\tau$ is as in the statement of the Proposition.
Hence $\beta p_{r} - \sigma_{2} \left( \beta p_{r} \right)=2b\sqrt{\tau}$ is an
algebraic integer and
$\left( \beta p_{r} - \sigma_{2} \left( \beta p_{r} \right) \right) \sqrt{\tau}$
is an algebraic integer in $\bK$. Similarly,
$(\gamma p_{r} - \sigma_{2} \left( \gamma p_{r} \right))\sqrt{\tau}$ is an
algebraic integer in $\bK$.

In this case, we put
$q_{r}' = \left( \gamma p_{r} - \sigma_{2} \left( \gamma p_{r} \right) \right)\sqrt{\tau}$
and $p_{r}' = \left( \beta p_{r} - \sigma_{2} \left( \beta p_{r} \right) \right)\sqrt{\tau}$
and observe that
$$
\left| \alpha q_{r}' - p_{r}' \right|
< l_{0} |\sqrt{\tau}| |\sigma_{2}(\beta)-\alpha \sigma_{2}(\gamma)| E^{-r}
$$
and
$$
\left| q_{r}' \right|
\leq k_{0} |\sqrt{\tau}| \left( |\gamma| + |\sigma_{2}(\gamma)| \right) Q^{r}.
$$

Note that in both these subcases, we obtain the same upper bound for
$\left| \alpha q_{r}' - p_{r}' \right|$ and for $\left| q_{r}' \right|$.

Here
$$
p_{r}'q_{r+1}'-p_{r+1}'q_{r}'
= \tau \left( \beta\sigma_{2}(\gamma) - \sigma_{2}(\beta)\gamma \right)
\left( \sigma_{2}(p_{r})p_{r+1} - p_{r}\sigma_{2}(p_{r+1}) \right),
$$
which we saw in Case~1 can only be zero if $\beta/\gamma \in \bK$ or
$p_{r}/p_{r+1} \in \bK$.

Therefore, we can apply Lemma~\ref{lem:approx} to find that
$\kappa = \log (Q)/\log (E)$ and
$$
c = 2k_{0} |\sqrt{\tau}| \left( |\gamma| + |\sigma_{2}(\gamma)| \right) Q
    \max \left\{ E, 2l_{0} |\sqrt{\tau}| |\sigma_{2}(\beta)-\alpha \sigma_{2}(\gamma)| E \right\}^{\kappa},
$$
concluding the proof of Case~2 and the Proposition.

\section{Proof of Theorem~\ref{thm:general-hypg}}

\subsection{Construction of approximations}

We construct the approximations under more general conditions. The point is not
to generalise for its own sake, but to illustrate the requirements and limitations
of our method of proof.

Let $\zeta_{k}$ be a $k$-th root of unity for some $k$. We apply
Lemma~\ref{lem:relation} with $z=\zeta_{k}\eta/\sigma(\eta)$. Multiplying both
sides of (\ref{eq:approx}) by $\sigma(\eta)^{r}$, we obtain
\begin{eqnarray*}
&   & \left( \zeta_{k} \eta/\sigma(\eta) \right)^{m/n}
      \left( \zeta_{k} \eta \right)^{r} X_{m,n,r} \left( \sigma(\eta)/(\zeta_{k}\eta) \right)
      - \sigma(\eta)^{r} X_{m,n,r} \left( \zeta_{k}\eta/\sigma(\eta) \right) \\
& = & \sigma(\eta)^{r} \left( \zeta_{k}\eta/\sigma(\eta)-1 \right)^{2r+1}
      R_{m,n,r} \left( \zeta_{k}\eta/\sigma(\eta) \right),
\end{eqnarray*}
which we can rewrite as
\begin{eqnarray*}
&   & \left( \zeta_{k}\eta/\sigma(\eta) \right)^{m/n}
      X_{m,n,r}^{*} \left( \sigma(\eta),\zeta_{k}\eta \right)
      - X_{m,n,r}^{*} \left( \zeta_{k}\eta, \sigma(\eta) \right) \\
& = & \sigma(\eta)^{r} \left( \zeta_{k}\eta/\sigma(\eta)-1 \right)^{2r+1}
      R_{m,n,r} \left( \zeta_{k}\eta/\sigma(\eta) \right).
\end{eqnarray*}

Observe that
\begin{eqnarray*}
X_{m,n,r}^{*} \left( \zeta_{k}\eta, \sigma(\eta) \right)
& = & g^{r} X_{m,n,r}^{*} \left( \frac{\zeta_{k}\eta}{g}, \frac{\sigma(\eta)}{g} \right) \\
& = & \left( g \frac{\sigma(\eta)}{g} \right)^{r}
      X_{m,n,r} \left( 1- d \frac{ (\sigma(\eta)-\zeta_{k}\eta)/g}{d\sigma(\eta)/g} \right).
\end{eqnarray*}

From Lemma~7.4(a) of \cite{Vout},
$$
\frac{D_{n,r}}{N_{d,n,r}} X_{m,n,r} \left( 1- d \frac{ (\sigma(\eta)-\zeta_{k}\eta)/g}{d\sigma(\eta)/g} \right)
\in \bZ \left[ \frac{ (\sigma(\eta)-\zeta_{k}\eta)/g}{d\sigma(\eta)/g} \right],
$$
and, as a consequence,
$$
\left( \frac{\sigma(\eta)}{g} \right)^{r}
\frac{D_{n,r}}{N_{d,n,r}} X_{m,n,r} \left( 1- d \frac{ (\sigma(\eta)-\zeta_{k}\eta)/g}{d\sigma(\eta)/g} \right)
$$
is an algebraic integer, since $(\sigma(\eta)-\zeta_{k}\eta)/(gd)$ is an algebraic
integer by the definition of $d$ in the statement of the Theorem. Hence
$$
p_{r} = \frac{h_{r}D_{n,r}}{g^{r} N_{d,n,r}}X_{m,n,r}^{*} \left( \zeta_{k}\eta, \sigma(\eta) \right)
$$
is an algebraic integer in $\bL$.

Similarly,
$$
q_{r} = \frac{h_{r}D_{n,r}}{g^{r} N_{d,n,r}}X_{m,n,r}^{*} \left( \sigma(\eta), \zeta_{k}\eta \right)
$$
is an algebraic integer in $\bL$.

Now we want $p_{r}$ and $q_{r}$, or at least numbers obtained from them, to be
algebraic conjugates. For this purpose, we must suppose that
$1/\zeta_{k}=\sigma(\zeta_{k})$ (note that this implies that $\zeta_{k} \in \bL$).

With this condition, and since $\sigma^{2}(\cdot)$ is the identity map, we have
\begin{eqnarray*}
\left( \zeta_{k} \right)^{r} \sigma \left( X_{m,n,r}^{*} \left( \zeta_{k} \eta, \sigma(\eta) \right) \right)
& = & \left( \zeta_{k} \right)^{r} \sigma \left( \sigma(\eta)^{r} X_{m,n,r} \left( \zeta_{k}\eta/\sigma(\eta) \right) \right) \\
& = & \left( \zeta_{k}\eta \right)^{r} X_{m,n,r} \left( \sigma \left( \zeta_{k}\eta/\sigma(\eta) \right) \right) \\
& = & \left( \zeta_{k}\eta \right)^{r} X_{m,n,r} \left( \sigma(\eta)/(\zeta_{k}\eta) \right) \\
& = & X_{m,n,r}^{*} \left( \sigma(\eta), \zeta_{k}\eta \right).
\end{eqnarray*}

Hence, $q_{r}=\zeta_{k}^{r}\sigma(p_{r})$ and so $q_{r}$ and $\sigma(\zeta_{k})^{r}p_{r}$
are algebraic conjugates over $\bK$. Letting, $k_{1}=k/(2,k)$, we have $p_{k_{1}r}$
and $\pm q_{k_{1}r}$ are algebraic conjugates for $k=1$, $2$, $3$, $4$ and $6$,
so we could put $p_{r}'=p_{k_{1}r}$ and $q_{r}'=q_{k_{1}r}$.

However here we restrict our attention to $k=1$ and observe that in this case
$p_{r}$ and $q_{r}$ are algebraic conjugates.

\subsection{Estimates}

From Lemmas~7.3(a) and 7.4(c) of \cite{Vout}, we have
\begin{eqnarray*}
\left| q_{r} \right|
& \leq & \frac{2h}{|g|^{r}}
         \frac{D_{n,r}}{N_{d,n,r}} \frac{\Gamma(1-m/n) r!}{\Gamma(r+1-m/n)}
         \max \left( \left| \sqrt{\eta}+\sqrt{\sigma(\eta)} \right|,
                     \left| \sqrt{\eta}-\sqrt{\sigma(\eta)} \right| \right)^{2r} \\
& \leq & 2h\cC_{n} \left( \frac{\cD_{n}}{|g|\cN_{d,n}} \right)^{r}
         \max \left( \left| \sqrt{\eta}+\sqrt{\sigma(\eta)} \right|,
                     \left| \sqrt{\eta}-\sqrt{\sigma(\eta)} \right| \right)^{2r}.
\end{eqnarray*}

From Lemma~7.2(a) of \cite{Vout},
\begin{eqnarray*}
&      & \left| \left( \sigma(\eta) \right)^{r} (\eta/\sigma(\eta)-1)^{2r+1} R_{m,n,r}(\eta/\sigma(\eta)) \right| \\
& \leq & 2.38 \left| 1- (\eta/\sigma(\eta))^{m/n} \right| \frac{n\Gamma(r+1+m/n)}{m\Gamma(m/n)r!} \\
&      & \times \min \left( \left| \sqrt{\eta}+\sqrt{\sigma(\eta)} \right|, \left| \sqrt{\eta}-\sqrt{\sigma(\eta)} \right| \right)^{2r}.
\end{eqnarray*}

Hence
\begin{eqnarray*}
&      & \left| q_{r} (\eta/\sigma(\eta))^{m/n} - p_{r} \right| \\
& \leq & 2.38h \frac{D_{n,r}}{|g|^{r}N_{d,n,r}} \left| 1- (\eta/\sigma(\eta))^{m/n} \right|
         \frac{n\Gamma(r+1+m/n)}{m\Gamma(m/n)r!} \\
&      & \times \min \left( \left| \sqrt{\eta}+\sqrt{\sigma(\eta)} \right|,
                            \left| \sqrt{\eta}-\sqrt{\sigma(\eta)} \right| \right)^{2r} \\
& \leq & \frac{2.38h}{|g|^{r}} \left| 1- (\eta/\sigma(\eta))^{m/n} \right|
         \cC_{n} \left( \frac{\cD_{n}}{\cN_{d,n}} \right)^{r} \\
&      & \times \min \left( \left| \sqrt{\eta}+\sqrt{\sigma(\eta)} \right|,
                            \left| \sqrt{\eta}-\sqrt{\sigma(\eta)} \right| \right)^{2r}.
\end{eqnarray*}

Therefore, in the notation of Proposition~\ref{prop:1}, we have
\begin{eqnarray*}
k_{0}  & = & 2h \cC_{n}, \\
l_{0}  & = & 2.38h \left| 1- (\eta/\sigma(\eta))^{m/n} \right| \cC_{n}, \\
E      & = & \left\{ \frac{\cD_{n}}{|g|\cN_{d,n}} \min \left( \left| \sqrt{\eta}-\sqrt{\sigma(\eta)} \right|^{2},
             \left| \sqrt{\eta}+\sqrt{\sigma(\eta)} \right|^{2} \right) \right\}^{-1}, \\
Q      & = & \frac{\cD_{n}}{|g|\cN_{d,n}} \max \left( \left| \sqrt{\eta}-\sqrt{\sigma(\eta)} \right|^{2},
             \left| \sqrt{\eta}+\sqrt{\sigma(\eta)} \right|^{2} \right).
\end{eqnarray*}

From Proposition~\ref{prop:1}, the expression for $\kappa$ in the Theorem
follows immediately, while, upon noting that our $\beta$, $\gamma$,
$\sigma(\beta)$ and $\sigma(\gamma)$ here are $\sigma_{2}(\beta)$,
$\sigma_{2}(\gamma)$, $\beta$ and $\gamma$ respectively in the notation of that
Proposition,
\begin{eqnarray*}
c & = & 2 |\sqrt{\tau}|\left( |\gamma| + |\sigma(\gamma)| \right) k_{0}Q
        \max \left\{ E, 2 |\sqrt{\tau}| \left( |\beta-\alpha \gamma| \right) l_{0}E \right\}^{\kappa} \\
  & < & 4h |\sqrt{\tau}| \left( |\gamma| + |\sigma(\gamma)| \right) \cC_{n}Q \\
  &   & \times \max \left\{ E, 5h |\sqrt{\tau}| \left| 1- (\eta/\sigma(\eta))^{m/n} \right|
                               |\beta-\alpha \gamma| \cC_{n} E \right\}^{\kappa}.
\end{eqnarray*}

\section{Proof of Theorem~\ref{thm:general-hypg-unitdisk}}

The proof of Theorem~\ref{thm:general-hypg-unitdisk} is the same as that of
Theorem~\ref{thm:general-hypg}, except that we use the upper bounds from
parts~(b) of Lemmas~7.2 and 7.3 of \cite{Vout}, rather than parts~(a).
Thus, we find that
\begin{eqnarray*}
k_{0}  & = & 2h \cC_{n}, \\
l_{0}  & = & h\left| 1- (\eta/\sigma(\eta))^{m/n} \right| \cC_{n}, \\
E      & = & \frac{4|g|\cN_{d,n}}{\cD_{n}} \frac{\left( |\eta| - |\sigma(\eta) - \eta | \right)}{|\sigma(\eta) - \eta |^{2}}, \\
Q      & = & \frac{2\cD_{n}}{|g|\cN_{d,n}} \left( \left| \eta \right| + \left| \sigma(\eta) \right| \right).
\end{eqnarray*}

So, from Proposition~\ref{prop:1}, $\kappa$ is as in the statement of the
Theorem and, again noting the change of notation mentioned at the end of the
proof of Theorem~\ref{thm:general-hypg},
\begin{eqnarray*}
c & = & 2 |\sqrt{\tau}| \left( |\gamma| + |\sigma(\gamma)| \right) k_{0}Q
        \max \left\{ E, 2 |\sqrt{\tau}| |\beta-\alpha \gamma| l_{0}E \right\}^{\kappa} \\
  & = & 4h |\sqrt{\tau}| \left( |\gamma| + |\sigma(\gamma)| \right) \cC_{n}Q \\
  &   & \times \max \left\{ E, 2h |\sqrt{\tau}| \left| 1- (\eta/\sigma(\eta))^{m/n} \right|
                               |\beta-\alpha \gamma| \cC_{n} E \right\}^{\kappa}.
\end{eqnarray*}

\section{Proof of Corollary~\ref{cor:cor-1}}

This Corollary follows from a direct application of Theorem~\ref{thm:general-hypg}.

We can write
\begin{equation}
\label{eq:ei}
\left( \sqrt{\eta} \pm \sqrt{\sigma(\eta)} \right)^{2}
= \eta + \sigma(\eta) \pm 2 \sqrt{\eta \sigma(\eta)}.
\end{equation}

The right-hand side of (\ref{eq:ei}) is
$u_{1} \pm \sqrt{u_{1}^{2}-u_{2}^{2}t}$ and $\sigma(\eta)-\eta
=-u_{2}\sqrt{t}$. Hence $d$ is as defined in the Corollary.

The analysis of $g_{1}$, $g_{2}$ and $g_{3}$ is identical to that in Section~11 of
\cite{Vout}.

As stated in the remark after Corollary~\ref{cor:cor-1}, $g_{4}$ and $g_{5}$
arise from the interplay of $d$ and $g$. Suppose that $d_{1}$ is the largest
positive rational integer such that $u_{2}\sqrt{t}/(d_{1}g_{1}\sqrt{g_{2}/g_{3}})$
is an algebraic integer. If there are multiplicative factors of the form
$\sqrt{d_{2}}$ in $u_{2}\sqrt{t}/(d_{1}g_{1}\sqrt{g_{2}/g_{3}})$, then by
multiplying $\eta$, and hence $u_{2}\sqrt{t}$, by $\sqrt{d_{2}}$, we can
increase $d_{1}$ by a factor of $d_{2}$. Under some circumstances, this
increases $\cN_{d,n}$ by a factor of $d_{2}$ while increasing
$u_{1} \pm \sqrt{u_{1}^{2}-u_{2}^{2}t}$ only by a factor of $\sqrt{d_{2}}$ for
a net reduction in the size of $\kappa$. We demonstrate here how $g_{4}$ and
$g_{5}$ capture these circumstances.

Consider the integer $u_{2}^{2}tg_{3}/(g_{1}^{2}g_{2})$ and let $d_{1}^{2}$ be
its largest square divisor. Suppose that $p$ is a prime divisor of their quotient.
That is, $p$ is a prime divisor of $\core(u_{2}^{2}tg_{3}/(g_{1}^{2}g_{2}))
=\core(tg_{3}/g_{2})=\core(tg_{2}g_{3})$. Note that
$$
d_{1}=\sqrt{u_{2}^{2}tg_{3}/(g_{1}^{2}g_{2})/\core(tg_{2}g_{3})}
=(u_{2}/g_{1})\sqrt{tg_{3}/g_{2}/\core(tg_{2}g_{3})}
$$.

First, if $p \nmid n$, then $\cN_{pd_{1},n}=\cN_{d_{1},n}$ from the definition
of $\cN_{d,n}$ in (\ref{eq:Ndn}) and there is no benefit.

Second, if $p|n$ and $p \nmid (n/\gcd(d_{1},n))$, then $\cN_{pd_{1},n}$ is at
most $\cN_{d_{1},n}p^{1/(p-1)}$ (again, from (\ref{eq:Ndn})). That is we gain
at most a factor of $p^{1/(p-1)}$, while increasing the size of
$u_{1} \pm \sqrt{u_{1}^{2}-u_{2}^{2}t}$ by a factor of $\sqrt{p}$ and hence
obtain no benefit for $p>2$.

Third, if $p|n$ and $p|(n/\gcd(d_{1},n))$, then we gain a factor of $p$, while
we increase the size of $u_{1} \pm \sqrt{u_{1}^{2}-u_{2}^{2}t}$ by a factor of
$\sqrt{p}$. The product of all such $p$ equals
$$
\gcd \left( \core(tg_{2}g_{3}), \frac{n}{\gcd((u_{2}/g_{1})\sqrt{tg_{3}/g_{2}/\core(tg_{2}g_{3})}, n)} \right),
$$
which is our $g_{4}$.

This covers all possible cases except $2|n$ and $2 \nmid (n/\gcd(d_{1},n))$.
If the power of $2$ dividing $d$ equals the power of $2$ dividing $n$, both are
positive and $2|\core(tg_{2}g_{3})$, then we increase $\cN_{d_{1},n}$ by a factor
of $2$, while we increase the size of $u_{1} \pm \sqrt{u_{1}^{2}-u_{2}^{2}t}$
by a factor of $\sqrt{2}$. Since
$u_{2}^{2}tg_{3}/(g_{1}^{2}g_{2})=d_{1}^{2}\core(tg_{2}g_{3})$, this condition
is equivalent to our condition in the definition of $g_{5}$.

Lastly, we must consider $h_{r}$ and $h$.

Since $g^{2} \in \bQ$, we can take $h_{r}=1$ for $r$ even. Since
$(g_{3}g_{4}g_{5}/g_{2})\core(g_{2}g_{3}g_{4}g_{5})$ is a perfect square, we
can take $h_{r}=\sqrt{\core(g_{2}g_{3}g_{4}g_{5})}$ for $r$ odd. Observe that
$g_{4}g_{5}|(2tg_{3}/g_{2})$, $g_{2}|t$ and $g_{3}|4$. Hence
$h_{r} \leq \sqrt{|2t|}$ for $r$ odd.

\subsection*{Acknowledgements}

The author is very grateful to the referee for their very attentive reading of
this paper and the corrections and improvements that they suggested. Their
diligence has led to a much better paper.
%
%The author also expresses his deep appreciation to the Judith and Lawrence Schulman
%Foundation for their generous support throughout the entire course of this work.
%

\end{document}